\documentclass{article}

\usepackage{latexsym}     % Pour certains symboles (\Box, ...)

\usepackage{BoxedEPS}
\SetepsfEPSFSpecial
\HideDisplacementBoxes
\newcommand{\rmax}{$\mathcal{R}_{\rm max}$}

\title{
{\bf Partial synchronicity and the
(max,+) semiring}}
\author{{Michael Mc Gettrick}}
\date{
}
\begin{document}
\maketitle
\section*{}
\noindent
{\bf Abstract. }
In this paper we illustrate how
non-stochastic (max,+) techniques
can be used to describe partial synchronization
in a Discrete Event Dynamical System. Our work
uses results from the spectral theory of dioids and
analyses (max,+) equations describing various
synchronization rules in a simple network.
The network in question is a transport network 
consisting of two routes joined at a single point, and
our Discrete Events are the departure times of transport
units along these routes. We calculate the maximum 
frequency of circulation of these units as a function of
the synchronization parameter. These functions allow us 
further to determine the waiting times on various routes,
and here we find critical parameters (dependent
on the fixed travel times on each route) which dictate the
overall behavoiur. We give explicit equations for these 
parameters and state the rules which enable optimal 
performance in the network (corresponding to minimum
waiting time).

\medskip
\section{Introduction}\label{1}

Considerable advances have been made in recent years in 
the algebraic description of Discrete Event Dynamical
Systems (DEDS). Chief among these is the approach using
the so-called (max,+) semiring \cite{gb}, with 
corresponding development of algorithms for spectral 
calculations \cite{num} and concrete applications
\cite{train}.

In this paper we illustrate how 
non-stochastic (max,+) techniques 
can be used to describe partial synchronization 
(in a manner to be defined below) in a DEDS. Our work
uses results from the spectral theory of dioids and 
analyses (max,+) equations describing various
synchronization rules in a simple network.

The paper is structured as follows. In Section \ref{2} we
review necessary mathematical properties of the 
(max,+) semiring. Section \ref{3} introduces the 
network on which our DEDS is based, and in Section \ref{4}
we describe the different (max,+) equations we
will solve. The resulting solutions lead, in Section \ref{5},
to critical phenomena in the network, and we calculate
which among the (max,+) equations give optimal results. 
A numerical example illustrates our 
results in Section \ref{6} we conclude in
Section \ref{7}. 

\section{(max,+) semiring}\label{2}

The (max,+) semiring (or dioid)
 \rmax\  is the set $\mathcal{R}\cup
\{ -\infty\}$ with the two operations $\oplus,\otimes$
defined by $a\oplus b={\rm max}(a,b)$ and 
$a\otimes b=a + b$. Note that both operations possess 
identity elements, $-\infty$ for $\oplus$
and 0 for $\otimes$, but while
$\otimes$ is invertible, $\oplus$ is idempotent:
$a\oplus a = a$. 
In the remainder of this paper we will often omit the
explicit multiplication symbols $\times, \otimes$, as is 
the convention.
Our convention will be that if an expression contains any
explicit symbol from  \rmax\
we will assume all such hidden symbols are also from \rmax :
$a+bc$ means $a+(b\times c)$, $a\oplus bc$ means
$a\oplus (b\otimes c)$, and $(a\oplus b)(c + d)$ means
$(a\oplus b)\otimes (c + d)$.

These operations can be extended to matrices 
$A,B\in \mathcal{R}_{\rm max}^{n\times n}$ by defining
$(A\oplus B)_{ij}=A_{ij}\oplus B_{ij}$ and
$(A\otimes B)_{ij}=\oplus_k(A_{ik}\otimes B_{kj}).$
As this may seem unusual to the novice reader, a small
example is in order.
\begin{eqnarray}\label{eqd}
&&
\left( \begin{array}{cc}
3 & -1\\
0 & 5\\
\end{array} \right)
\otimes
\left( \begin{array}{c}
4\\
0\\
\end{array} \right)\nonumber\\
&&=
\left( \begin{array}{c}
(3\otimes 4)\oplus (-1\otimes 0)\\
(0\otimes 4)\oplus (5\otimes 0)\\
\end{array} \right)
=
\left( \begin{array}{c}
7\\
5\\
\end{array} \right).
\end{eqnarray}
We now look at solutions to
the (max,+) spectral problem
$A\otimes X = \lambda\otimes X$.\\\\
{\bf Theorem 1.}\cite{gun}\cite{maslov}
{\it Let $A\in\mathcal{R}_{\rm max}^{n\times n}$.
There exists a maximal eigenvalue }
\begin{eqnarray}\label{eqe}
\lambda =
\bigoplus_{k=1}^n
\bigoplus_{i_1,i_2,\dots,i_k}
(A_{i_1i_2}\otimes A_{i_2i_3}\otimes\dots
\otimes A_{i_ki_1})/k.
\end{eqnarray}
If we extend our notation to
$a^n
\stackrel{\rm def}{=}
a\otimes a\otimes\dots\otimes a = n\times a$ and
${\rm Tr}(A)=\bigoplus_{k=1}^n A_{kk}$ 
equation (\ref{eqe}) reads
\begin{eqnarray}\label{eqf}
\lambda=
\bigoplus_{k=1}^n
({\rm Tr}(A^k))^{1/k}.
\end{eqnarray}
In the standard way, we associate with $A$ a 
weighted digraph (precedence graph)
$\mathcal{G}(A)$ with nodes $\mathcal{N}=\{1,2,\dots,n\}$
and edges $\mathcal{E}=\{(i,j)|A_{ji}\neq -\infty\}$.
The weight of edge $(i,j)$ is simply $A_{ji}$. Denoting
by $i\to j$ the existence of a path from $i$ to $j$ we 
define $A$ to be irreducible (and $\mathcal{G}(A)$ to be 
strongly connected) iff $i\to j\quad\forall\ 1\leq i,j\leq n$.\\\\
{\bf Theorem 2.}\cite{gb}
{\it
If $A\in\mathcal{R}_{\rm max}^{n\times n}$ 
is irreducible, then $\lambda$ as given by 
{\rm (\ref{eqe})} is unique.}\\\\
Note that from a graph-theoretic point of view,
$\lambda$ is merely the maximum cycle mean of
$\mathcal{G}(A)$.

Our DEDS will be described by equations in
\rmax: Specifically, denoting the time of the
$k$th occurence of event
$i$ by $x_i(k)$, we are interested in equations of the form

\begin{eqnarray}\label{eqa}
x_i(k)=f(x_j(k-l))
\end{eqnarray}
where we have a finite family of $1\leq i,j\leq n$ events
and $1\leq l\leq k-1$. For linear equations in \rmax, 
(\ref{eqa}) becomes

\begin{eqnarray}\label{eqb}
x_i(k)=\bigoplus_{j,l}(f^l_{ij}\otimes x_j(k-l))
\end{eqnarray}
and we will show in Section \ref{4} how this can be 
written as a matrix equation with unitary retard

\begin{eqnarray}\label{eqc}
x_i(k)=\bigoplus_j(A_{ij}\otimes x_j(k-1)).
\end{eqnarray}
We are interested in the asymptotic properties of this
system as $k\to\infty$ so we study $\lim_{k\to\infty}A^k$.\\\\
{\bf Theorem 3.}\cite{steph}\cite{mplus}\cite{cun}
{\it If $A\in\mathcal{R}_{\rm max}^{n\times n}$ 
is irreducible, there exists integers $M, c(A)$
such that
\begin{eqnarray}\label{eqg}
A^{k\otimes c(A)} = (\lambda)^{c(A)}\otimes 
A^k\quad\forall\ k\geq M.
\end{eqnarray}}

We call $c(A)$ the cyclicity of $A$ (see \cite{gb} for
details of how to calculate this integer).

\section{Model}\label{3}

We consider an elementary model of two tourbuses with partial 
synchronization. Each tourbus goes on a fixed circular route $R_1$ 
and $R_2$, stopping every so often to drop off and pick up passengers.
The tour buses meet at one location, a downtown station S, where 
passengers can pass from one bus to the other one. The physical network
consists of two loops meeting at S, as shown in Figure 1.

\begin{figure}[h]
%[hbtp]
$$
\BoxedEPSF{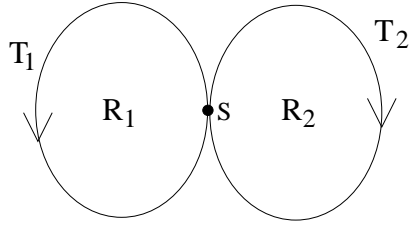 scaled 600}
$$
\caption{{\tt Physical} network of routes $R_1, R_2$ linked
at station S.}
\end{figure}

We assume fixed travel times $T_1$ and $T_2$ 
for $R_1$ and $R_2$ respectively, and without
loss of generality let $T_1 \le T_2$. 
The simplest model is to describe this
system using the (max,+) equations
$X(k+1) = A\otimes X(k)$, with

\begin{eqnarray}\label{eq1}
A= \left( \begin{array}{cc}
a_{11} & a_{12} \\
a_{21} & a_{22} \\
\end{array}
\right),\
X(k)=
\left( \begin{array}{c}
x_1(k)\\
x_2(k)\\
\end{array}
\right)
\end{eqnarray}
and $x_i(k)$ is the departure time of bus $i$
from station $S$ along route $R_i$.
We consider the following extreme cases.\\\\
%\begin{itemize}
%\item
{\bf CASE M1 (no synchronization):}
\par
We have
\begin{eqnarray}\label{eqh}
A= \left( \begin{array}{cc}
T_1 & -\infty\\
-\infty &T_2\\
\end{array}
\right)
\end{eqnarray}
so that the graph $\mathcal{G}(A)$ is disconnected into
two parts each having its own cycle mean 
$\lambda_1 = T_1$ and $\lambda_2 = T_2$ as shown in
Figure 2 (the nodes here correspond to {\bf departure times } of
buses).
\begin{figure}[h]
%[hbtp]
$$
\BoxedEPSF{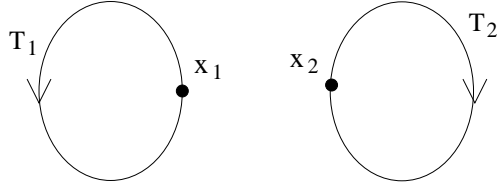 scaled 500}
$$
\caption{Precedence Graph. One bus each on $R_1$ and $R_2$.
Not Synchronized.}
\end{figure}\\\\
%\item
{\bf CASE M2 (complete synchronization):}
\par
We have
\begin{eqnarray}\label{eqh}
A= \left( \begin{array}{cc}
T_1 & T_2\\
T_1 &T_2\\
\end{array}
\right)
\end{eqnarray}
so that the graph $\mathcal{G}(A)$ is strongly 
connected.
There are 3 circuits in this graph with cycle means $T_1$,
$T_2$, and ${T_1 + T_2 \over 2}$, so the maximum cycle mean gives
an eigenvalue of
$\lambda = T_2$ (see Figure 3). 

%\end{itemize}

\begin{figure}[h]
%[hbtp]
$$
\BoxedEPSF{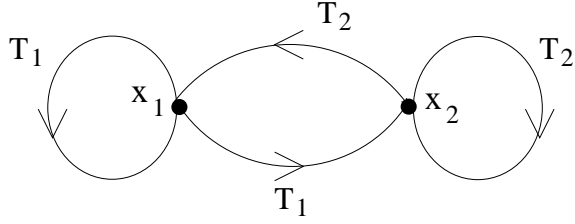 scaled 600}
$$
\caption{Precedence Graph. One bus each on $R_1$ and $R_2$.
Completely Synchronized.}
\end{figure}

\section{Partial Synchronicity}\label{4}

For the case of complete synchronization in section
\ref{3}, it is clear that the large value of $\lambda
= T_2$ slows down the overall system. We further suppose in
this section that we cannot alter the travel time $T_2$, but
that we have a large amount of additional stock (buses) which
we can add to the network.

We add to route $R_2$ a large number $m$ of buses such that 
the interval between departure times on this route becomes a
very small number $\delta$ (which we can think of as $T_2/m$ if
required). 
Denote by $y_1(k), y_2(k),\dots,y_m(k)$ the departure times of the
$m$ new buses, which operate according to the rules
\begin{eqnarray}\label{eq24}
y_i(k)=&\delta\otimes  y_{i-1}(k-1)\nonumber\\
y_0(k)=&\delta\otimes y_m(k-1).
\end{eqnarray}
For this route we
do not produce a timetable, as we assume to a first approximation
there is virtually no waiting for a bus.

For this system, it would seem there should be no interdependence
between $x_1(k)$ and $x_2(k)$. However, what of a passenger who 
wants to make a non-stop round trip $R_1\to R_2\to R_1$ -- would 
some synchronization rule benefit such a passenger? 
To investigate further we set $T_2=nT_1+r$ and look
at the \rmax\ 
equation
$x_1(k)=f(x_1(k-l))$ for different values of $l$ and 
linear $f$.\\\\
%\begin{itemize}
%\item
{\bf CASE P1:}
\begin{eqnarray}\label{eq23}
x_1(k)=T_1x_1(k-1)\oplus T_1T_2x_1(k-(l+1))
\end{eqnarray}
Using $n$ artificial states (i.e. states
that do not in fact correspond to the departure of any
physical bus) $z_1, z_2,\dots , z_l$, we can rewrite
(\ref{eq23}) as
\begin{eqnarray}\label{eq24}
x_1(k)&=&T_1x_1(k-1) \oplus rT_1z_l(k-1)\nonumber\\
z_i(k)&=&T_1^{n/l}\otimes z_{i-1}(k-1), \bigskip i=1,\dots,l\\
z_0(k)&\stackrel{\rm def}{=}
& x_1(k).\nonumber
\end{eqnarray}
At the expense of increasing the dimension, our 
system of equations now has unitary retard as in
(\ref{eqc}) with column vector
\begin{eqnarray}\label{eq24}
X(k)=
\left( \begin{array}{c}
x_1(k)\\
z_l(k)\\
z_{l-1}(k)\\
\vdots\\
z_1(k)\\
\end{array} \right)
\end{eqnarray}
and matrix
\begin{eqnarray}\label{eq25}
A=
\left( \begin{array}{cccccc}
T_1 & T_1 + r & -\infty & -\infty &\cdots &-\infty\\
-\infty & -\infty & nT_1/l & -\infty &\cdots &-\infty\\
-\infty & -\infty & -\infty & nT_1/l &\cdots &-\infty\\
\vdots & \vdots & \vdots & \vdots & \ddots & \vdots \\
nT_1/l & -\infty  & -\infty  & -\infty & \cdots  & -\infty \\
\end{array} \right).\nonumber
\end{eqnarray}
The full precedence graph  $\mathcal{G}(A)$
for the system $x_1, z_1,\dots
,z_l, x_2, y_1,\dots , y_m$ is shown in Figure 4.
\begin{figure}[h]
%[hbtp]
$$
\BoxedEPSF{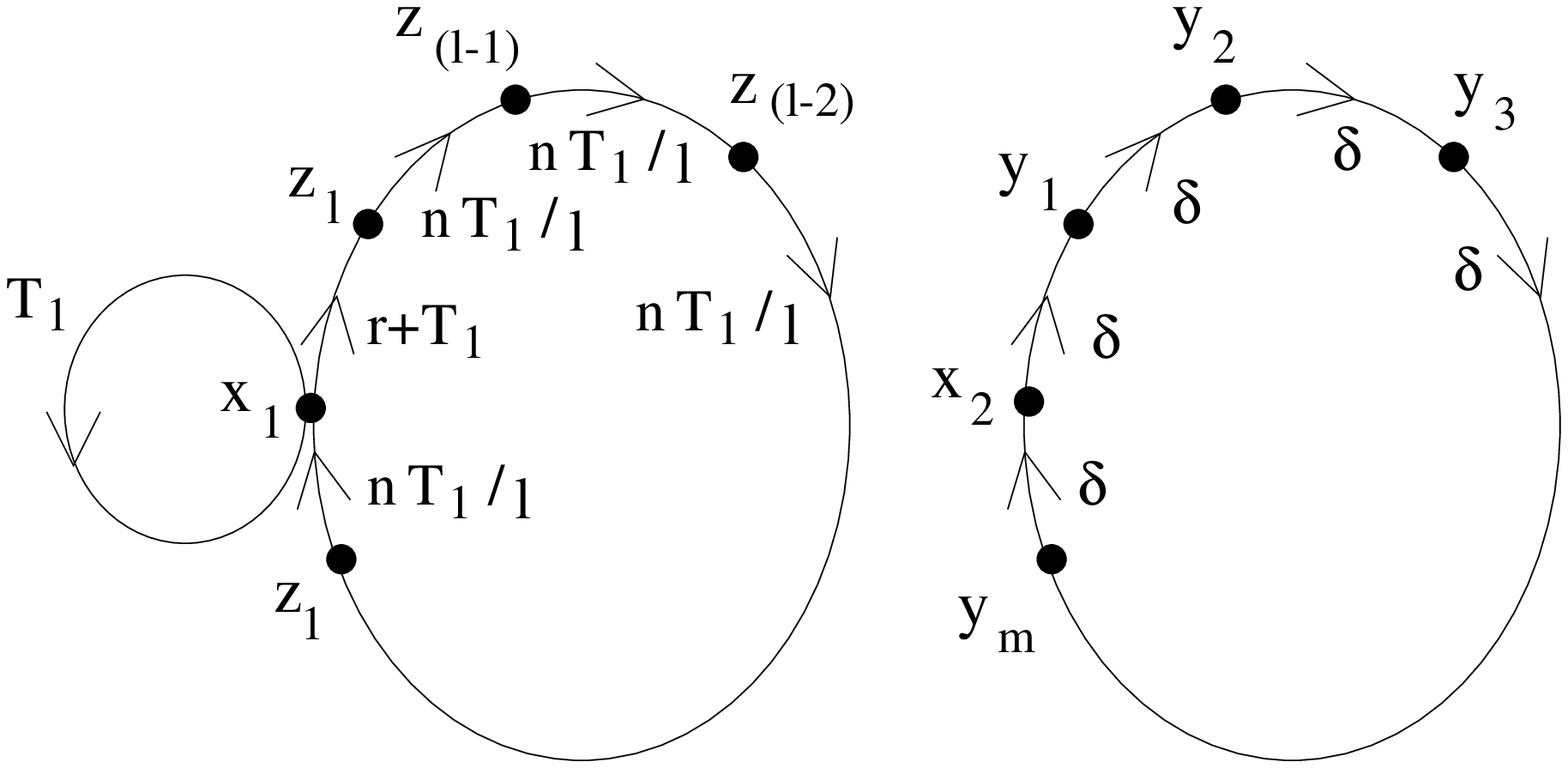 scaled 430}
$$
\caption{Precedence Graph for CASE P1, 
equation (\ref{eq23})}
\end{figure}
The disconnected parts of the graph $\mathcal{G}(A)$ 
have maximum cycle means $\lambda_2=\delta$ and
\begin{eqnarray}\label{eq26}
\lambda_1=
\left\{
\begin{array}{cc}
T_1 & \hbox{if } l>n\\
\frac{(n+1)T_1+r}{l+1} & \hbox{otherwise}
\end{array}
\right.
\end{eqnarray}
since $T_1\oplus r=T_1$ (by definition of
$r$) and, for $l>n$ (equivalently $l=l\oplus (n+1)$),
\begin{eqnarray}\label{eq32}
&\frac{(n+1)T_1+r}{l+1} <
\frac{(n+1)T_1+T_1}{l+1}  = \frac{(n+2)T_1}{l+1}\leq T_1. &
\end{eqnarray}
Note that $ l>n $ effectively corresponds to no
synchronization, while for the choice $ l\leq n $,
$(n-l+1)$ out of every $(n+1)$ buses are synchronized.
The minimum possible $\lambda_1$ is $\lambda_1=T_1$
given by $l>n$.\\\\
%\item
{\bf CASE P2:}
\begin{eqnarray}\label{eq27}
x_1(k)=T_1x_1(k-1)\oplus rT_1^{l+1}x_1(k-(l+1))
\end{eqnarray}
As in the previous case P1, we can rewrite (\ref{eq27})
as 
\begin{eqnarray}\label{eq28}
x_1(k)&=&T_1x_1(k-1) \oplus rT_1z_l(k-1)\nonumber\\
z_i(k)&=&T_1\otimes z_{i-1}(k-1), \bigskip i=1,\dots,l\\
z_0(k)&\stackrel{\rm def}{=}
& x_1(k).\nonumber
\end{eqnarray}
and the system is described by $X(k+1)=A\otimes X(k)$
with $X(k)$ as in (\ref{eq24}) and
\begin{eqnarray}\label{eq29}
A=
\left( \begin{array}{cccccc}
T_1 & T_1 + r & -\infty & -\infty &\cdots &-\infty\\
-\infty & -\infty & T_1 & -\infty &\cdots &-\infty\\
-\infty & -\infty & -\infty & T_1 &\cdots &-\infty\\
\vdots & \vdots & \vdots & \vdots & \ddots & \vdots \\
T_1 & -\infty  & -\infty  & -\infty & \cdots  & -\infty \\
\end{array} \right).\nonumber
\end{eqnarray}
From the precedence graph shown in Figure 5 we have
$\lambda_2=\delta$ and
\begin{eqnarray}\label{eq30}
\lambda_1=
\frac{(l+1)T_1+r}{l+1}
\end{eqnarray}
since obviously
\begin{eqnarray}\label{eq33}
\frac{(l+1)T_1+r}{l+1}\leq T_1 
\end{eqnarray}
for all values of $l$.
\begin{figure}[h]
%[hbtp]
$$
\BoxedEPSF{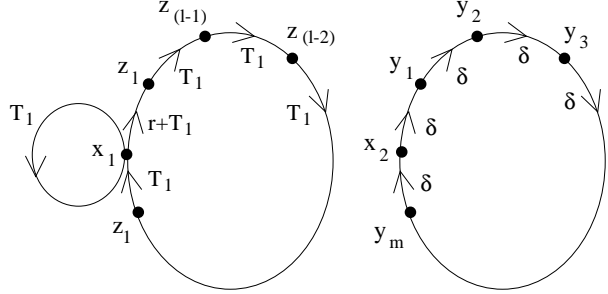 scaled 430}
$$
\caption{Precedence Graph for CASE P2, equation
(\ref{eq27})}
\end{figure}
This system corresponds to one out of every $(l+1)$ buses
being synchronized. The limiting case is
\begin{eqnarray}\label{eq31}
\lim_{l\to\infty}  \frac{(l+1)T_1+r}{l+1} = T_1
\end{eqnarray}
so the minimum value of $\lambda_1$ occurs at 
$l=\infty$.

\section{Optimal Solutions}\label{5}

In both cases P1 and P2 of Section \ref{4} the minimum 
value of $\lambda_1$ is $T_1$: In P1 it is given by any 
choice $l>n$ while in P2 it is given by the extreme value
$l=\infty$. These values both correspond to no synchronization,
and this is the result we would expect, that the system runs
faster if the buses do not wait for each other.

If our goal is to choose the synchronization rule which 
maximizes the speed of the system, the problem is solved.
In this Section we look at the separate problem of minimizing
the waiting time of passengers that make a circular journey over
all or part of the network (by circular we mean a journey that
begins and ends at the same point).
We itemize such journeys $J_i$ as follows:
\begin{eqnarray}\label{eq34}
J_1 &:& R_1\to R_1\nonumber\\
J_2 &:& R_1\to R_2\to R_1\nonumber\\
J_3 &:& R_2\to R_1\to R_2\nonumber\\
J_4 &:& R_2\to R_2\nonumber
\end{eqnarray}
For cases P1, P2 we look at regions $l\leq n$ and $l>n$ to 
give synchronization rules $S_i$ as follows:
\begin{eqnarray}\label{eq35}
S_1 &:& \hbox{Equation (\ref{eq23})}, l\leq n\nonumber\\
S_2 &:& \hbox{Equation (\ref{eq23})}, l>n\nonumber\\
S_3 &:& \hbox{Equation (\ref{eq27})}, l\leq n\nonumber\\
S_4 &:& \hbox{Equation (\ref{eq27})}, l>n\nonumber
\end{eqnarray}
For $i, j \in \{1,2,3,4\}$ we pair
each journey $J_i$ with  rules $S_j$
to  give 16 different models. Let $w_{ij}$ be the
total waiting time on journey $J_i$ with rule $S_j$.
These waiting times can be broken into two components,
the waiting time before boarding the bus  $w^b_{ij}$  and
the waiting time in transit $w^t_{ij}$, which can only occur 
for journeys $J_2$ and $J_3$. From our calculation of 
$\lambda_1$ and $\lambda_2$ in Section \ref{4} we determine 
$w^b_{ij}, w^t_{ij}$ and $w_{ij} = w^b_{ij} + w^t_{ij}$ to
be the values given by equations (\ref{eq35})--(\ref{eq37})
in the Appendix.

We now make some observations on the solutions given by
equation (\ref{eq37}). Our goal is to minimize the waiting times,
or more specifically the average waiting time per passenger
as a function of the synchronization parameter $l$:
If $n_i$ people take journey $J_i$ and $N=\sum_i n_i$
then we want to find the
minimum of
\begin{eqnarray}\label{eq38}
\overline{W_i}\stackrel{\rm def}{=}
\frac{\sum_j n_jw_{ji}}{N}
\end{eqnarray}
Note that
\begin{itemize}
\item
$w_{i1}$ is minimized by choosing the largest possible
value of $l$, i.e. $l=n$.

\item
$w_{i2}=\lim_{l\to\infty}w_{i4}$

\item
similarly $w_{i3}$ is minimized by choosing 
$l=n$.

\end{itemize}
Finally examine $w_{i4}$.
$w_{14}$ and $w_{34}$ are both clearly minimized by $l=\infty$.
To minimize $w_{24}$ we view it as a continuous function
of $l$ and solve
\begin{eqnarray}\label{eq39}
\frac{dw_{24}}{dl}=
\frac{d}{dl}\frac{(3l-2n+1)T_1 + (2n-2l+3)r}{2(l+1)}
=0.
\end{eqnarray}
$l$ drops out of equation (\ref{eq39}), meaning the existence of
an extremum does not depend on $l$ but rather on the solution of
$(2n+2)T_1=(2n+5)r$.
Thus we have a critical point 
\begin{eqnarray}\label{eq40}
r_c=\left(\frac{2n+2}{2n+5}\right)T_1.
\end{eqnarray}
The slope of the function $w_{24}(l)$ can now
be written in terms of the critical remainder $r_c$ as
\begin{eqnarray}\label{eq41}
\frac{dw_{24}}{dl}=
\frac{(2n+5)(r_c-r)}{(l+1)^2}.
\end{eqnarray}
Hence $w_{24}(l)$ is
\begin{itemize}
\item
strictly increasing for $r<r_c$: minimum is obtained by
choosing lowest value of $l$, i.e. $ l=n+1$
\item
strictly decreasing for $r>r_c$: minimum is obtained by
choosing highest value of $l$, i.e. $ l=\infty$
\item
constant for $r=r_c$: choose any value of $l$
\end{itemize}
In equation (\ref{eq42}) we rewrite $W$ with the optimal
values of $l$ chosen as indicated above: for $w_{i4}$ we 
insert in the matrix $W^m (\equiv W^{\min})$ the values at $l=n+1$, 
since the values at $l=\infty$ correspond to column 2,
$w_{i2}$. 

Now that we have chosen $l$ we should choose the
synchronization rule $S_j$ that minimizes 
$\overline{W_i}$. In rows 1, 3 and 4 of $W^m$ it is
clear that the minimum values are in column 2 corresponding
to $S_2$. Let us look at row 2, $w^m_{2j}$. In Figure 6 we plot
$w^m_{21}, w^m_{22} $ and $ w^m_{24}$ as functions of $r$.
\begin{figure}[h]
%[hbtp]
$$
\BoxedEPSF{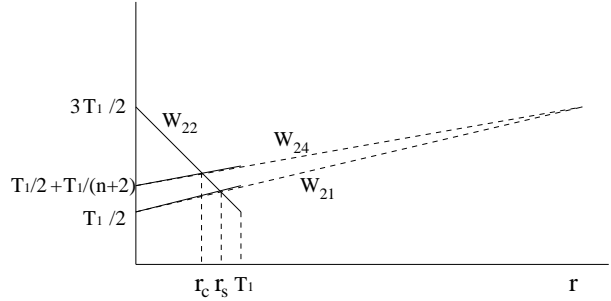 scaled 550}
$$
\caption{ $w^m_{21}, w^m_{22} $ and $ w^m_{24}$
plotted against $0\leq r\le T_1$ (see equation
(\ref{eq42})): The critical
points are $r_c=(2n+2)T_1/(2n+5)$ and 
$r_s=(2n+2)T_1/(2n+3)$.
}
\end{figure}
\begin{eqnarray}\label{eq44}
r_s=\left(\frac{2n+2}{2n+3}\right)T_1
\end{eqnarray}
along with equation (\ref{eq40})
defines three different regimes as follows:
\begin{itemize}
\item
$r>r_s: w^m_{22}<w^m_{21}<w^m_{24}$

\item
$r_c<r<r_s: w^m_{21}<w^m_{22}<w^m_{24}$

\item
$r<r_c: w^m_{21}<w^m_{24}<w^m_{22}$

\end{itemize}
Further we note the importance of the critical value
$r_s$ is that for $r>r_s$ rule $S_2$ gives shortest waiting
times for all journeys $J_i$, and in this regime 
$\overline{W_2}$ is minimal irrespective of the values
$n_i$ (see equation (\ref{eq38})). For $r<r_s$  the values
of $n_i$ dictate which $\overline{W_i}$ is minimal:
In particular, if
\begin{eqnarray}\label{eq43}
n_2>\frac{r}{(2n+3)(r_s-r)}(3n_1+n_3)
\end{eqnarray}
then $\overline{W_1}$ is minimal. We now have a rule that governs
optimum performance in our network: If $r>r_s$ implement 
rule $S_2$, otherwise monitor passenger numbers $n_i$ and 
if at any stage equation (\ref{eq43}) holds, implement
rule $S_1$.

Observe further from equation (\ref{eq44})
that $4T_1/5\leq r_s < T_1$, since $1\leq n < \infty$. This 
means that even for a random choice of $T_1$ and $T_2$, the
chances that $r$ is less than $r_s$ are greater than 80\%:
\begin{eqnarray}\label{eq45}
P(r<r_s | \hbox{ arbitrary } T_1, T_2) \geq 0.8
\end{eqnarray}

\section{Example}\label{6}

Let $T_1 = 3$ units of time and $n=2$ so that 
$r_s=18/7$ as given by equation (\ref{eq44}).

To examine firstly the regime $r>r_s$, let 
$r=8/3$, given by $T_2=26/3$ units of time. Under rule
$S_2$ we calculate $w^m_{22}=33/18$ with actual departure
times $x_1(k)=0, 3, 6, 9, 12, 15,\dots $ Under rule $S_3$
we have $w^m_{23}=35/18$ with departure times 
$x_1(k)=0, 3, 6, 35/3, 44/3, 53/3,\dots $ Thus rule $S_2$
gives optimal performance on $J_2$.

Secondly, let $r=7/3$ given by $T_2=25/3$, so that we are in
the regime $r<r_s$. Rule $S_2$ gives $w^m_{22}=39/18$ with
departure times $x_1(k)=$ 0, 3, 6, 9, 12, 15,$\dots $ Under rule $S_3$
we have $w^m_{23}=34/18$ with departure times 
$ x_1(k)= $ 0, 3, 6, 34/3, 43/3, 52/3,$\dots $ Here optimal performance
on $J_2$
is given by rule $S_3$.

Note in both regimes of this example the departure times $x_1(k)$
under $S_2$ are identical, the reason for the differences in 
$w^m_{22}$ lying in the transfer waiting time.

\section{Conclusion}\label{7}

In this paper we have solved the (max,+) equations for a 
specific model with two routes joining at one point. We have
shown how partial synchronicity can be described by these 
equations by introducing a synchronization parameter $l$.
Using results from (max,+) spectral theory we have calculated
the maximum frequency (minimum cycle mean) of buses on such 
routes. These numbers allow us to calculate waiting times on 
various routes: We observe new critical behaviour which 
depends on the relation between $T_1$ and $T_2$, the (fixed)
travel times on each route. Finally we state the rules for 
optimal performance in this network.

Further work is envisaged in looking at more general (max,+)
equations, in particular a system with $m_1$ buses on $R_1$ and
$m_2$ on $R_2$. Another direction in which this work could possibly
be extended is to have further routes $R_3, R_4,\dots$

\section*{Appendix}
\begin{eqnarray}\label{eq35}
W^b=
\left( \begin{array}{cccc}
\vspace{0.2cm}
\frac{(n+1)T_1+r}{2(l+1)} & \frac{T_1}{2} &
\frac{(l+1)T_1+r}{2(l+1)} & \frac{(l+1)T_1+r}{2(l+1)}\\
\vspace{0.2cm}
\frac{(n+1)T_1+r}{2(l+1)} & \frac{T_1}{2} &
\frac{(l+1)T_1+r}{2(l+1)} & \frac{(l+1)T_1+r}{2(l+1)}\\
\vspace{0.2cm}
0 & 0 & 0 & 0\\
0 & 0 & 0 & 0\\
\end{array} \right)
\end{eqnarray}
\begin{eqnarray}\label{eq36}
W^t=
\left( \begin{array}{cccc}
\vspace{0.2cm}
\frac{(n-l)T_1+r}{l+1} & 0 &
\frac{r}{l+1} & \frac{r}{l+1}\\
\vspace{0.2cm}
0 & T_1-r & 
\frac{(n-l)r}{(l+1)} & \frac{T_1 +(T_1-r)[l-(n+1)]}{(l+1)}\\
\vspace{0.2cm}
\frac{(n+1)T_1+r}{2(l+1)} & \frac{T_1}{2} &
\frac{(l+1)T_1+r}{2(l+1)} & \frac{(l+1)T_1+r}{2(l+1)}\\
0 & 0 & 0 & 0\\
\end{array} \right)
\end{eqnarray}
\begin{eqnarray}\label{eq37}
W=W^b+W^t=
\left( \begin{array}{cccc}
\vspace{0.2cm}
\frac{(3n-2l+1)T_1+3r}{2(l+1)} & \frac{T_1}{2} &
\frac{(l+1)T_1+3r}{2(l+1)} & \frac{(l+1)T_1+3r}{2(l+1)}\\
\vspace{0.2cm}
\frac{(n+1)T_1+r}{2(l+1)} & \frac{3T_1}{2}-r & 
\frac{(l+1)T_1+(2n-2l+1)r}{2(l+1)} & 
\frac{(3l-2n+1)T_1 + (2n-2l+3)r}{2(l+1)}\\
\vspace{0.2cm}
\frac{(n+1)T_1+r}{2(l+1)} & \frac{T_1}{2} &
\frac{(l+1)T_1+r}{2(l+1)} & \frac{(l+1)T_1+r}{2(l+1)}\\
0 & 0 & 0 & 0\\
\end{array} \right)
\end{eqnarray}
\begin{eqnarray}\label{eq42}
W^m (\equiv W^{\min}) =
\left( \begin{array}{cccc}
\vspace{0.2cm}
\frac{T_1}{2}+\frac{3r}{2(n+1)} & \frac{T_1}{2} &
\frac{T_1}{2}+\frac{3r}{2(n+1)} & \frac{T_1}{2} +\frac{r}{2(n+2)}\\
\vspace{0.2cm}
\frac{T_1}{2}+\frac{r}{2(n+1)} & \frac{3T_1}{2}-r &
\frac{T_1}{2}+\frac{r}{2(n+1)} &
\frac{T_1}{2}+\frac{2T_1+r}{2(n+2)}\\
\vspace{0.2cm}
\frac{T_1}{2}+\frac{r}{2(n+1)} & \frac{T_1}{2} &
\frac{T_1}{2}+\frac{r}{2(n+1)} & \frac{T_1}{2}+\frac{r}{2(n+2)}\\
0 & 0 & 0 & 0\\
\end{array} \right)
\end{eqnarray}

\noindent
Michael Mc Gettrick\\
Department of Information Technology\\
National University of Ireland\\
Galway\\
Ireland\\
e-mail: michael.mcgettrick@nuigalway.ie

\end{document}